\documentclass[12pt]{amsart}
\usepackage{eufrak}
\usepackage[total={6.5in,8.75in}, top=1.2in, left=1.0in, includefoot]{geometry}
\usepackage{ amssymb }

\theoremstyle{definition}

\theoremstyle{remark}

\numberwithin{equation}{section}

\begin{document}

\title{Specialization of monodromy group and $\ell$-independence}
\author{Chun Yin Hui}

\begin{abstract}

Let $E$ be an abelian scheme over a geometrically connected variety $X$ defined over $k$, a finitely generated field over $\mathbb{Q}$. Let $\eta$ be the generic point of $X$ and $x\in X$ a closed point. If $\mathfrak{g}_l$ and $(\mathfrak{g}_l)_x$ are the Lie algebras of the $l$-adic Galois representations for abelian varieties $E_{\eta}$ and $E_x$, then $(\mathfrak{g}_l)_x$ is embedded in $\mathfrak{g}_l$ by specialization. We prove that the set $\{x\in X$ closed point $|~ (\mathfrak{g}_l)_x\subsetneq \mathfrak{g}_l\}$ is independent of $l$ and confirm Conjecture 5.5 in [2].

\end{abstract}

\maketitle

\textbf{$\mathsection0$. Introduction}\\

Let $E$ be an abelian scheme of relative dimension $n$ over a geometrically connected variety $X$ defined over $k$, a finitely generated field over $\mathbb{Q}$. If $K$ is the function field of $X$ and $\eta$ is the generic point of $X$, then $A:=E_{\eta}$ is an abelian variety of dimension $n$ defined over $K$. The structure morphism $X\rightarrow \mathrm{Spec}(k)$ induces at the level of $\acute{e}tale$ fundamental groups a short exact sequence of profinite groups:
\begin{equation}
1\rightarrow \pi_1(X_{\overline{k}})\rightarrow \pi_1(X)\rightarrow \Gamma_k:=\mathrm{Gal}(\overline{k}/k)\rightarrow 1.
\end{equation}
 Any closed point $x:\mathrm{Spec}(\textbf{k}(x))\rightarrow X$ induces a splitting $x:\Gamma_{\textbf{k}(x)}\rightarrow \pi_1(X_{\textbf{k}(x)})$ of equation $(0.1)$ for $\pi_1(X_{\textbf{k}(x)})$.\\

Let $\Gamma_K=\mathrm{Gal}(\overline{K}/K)$ the absolute Galois group of $K$. For each prime number $l$, we have the Galois representation $\rho_{l}:\Gamma_K\rightarrow \mathrm{GL}(T_l(A))$ where $T_l(A)$ is the $l$-adic Tate module of $A$. This representation is unramified over $X$ and factors through $\rho_l:\pi_1(X)\rightarrow \mathrm{GL}(T_l(A))$ (still denote the map by $\rho_l$ for simplicity). The image of $\rho_{l}$ is a compact $l$-adic Lie subgroup of $\mathrm{GL}(T_l(A))\cong \mathrm{GL}_{2n}(\mathbb{Z}_l)$. Any closed point $x:\mathrm{Spec}(\textbf{k}(x))\rightarrow X$ induces an $l$-adic Galois representation by restricting $\rho_l$ to $x(\Gamma_{\textbf{k}(x)})$. This representation is isomorphic to the Galois representation of $\Gamma_{\textbf{k}(x)}$ on the $l$-adic Tate module of $E_x$, the abelian variety specialized at $x$.\\

For simplicity, write $G_l:=\rho_l(\pi_1(X))$, $\mathfrak{g}_l:=\mathrm{Lie}(G_l)$, $(G_l)_x:=\rho_l(x(\Gamma_{\textbf{k}(x)}))$ and $(\mathfrak{g}_l)_x:=\mathrm{Lie}((G_l)_x)$. We have $(\mathfrak{g}_l)_x\subset \mathfrak{g}_l$. We set $X^{cl}$ the set of closed points of $X$ and define the exceptional set
\begin{equation*}
X_{\rho_{E,l}}:=\{x\in X^{cl}| (\mathfrak{g}_l)_x\subsetneq \mathfrak{g}_l\}.
\end{equation*}
The main result (Theorem 1.4) of this note is that the exceptional set $X_{\rho_{E,l}}$ is independent of $l$. Conjecture 5.5 in [Cadoret \& Tamagawa 2] is then a direct application of our theorem.\\

\textbf{$\mathsection1$. $l$-independence of $X_{\rho_{E,l}}$}\\

\textbf{Theorem 1.1.} (Serre [5 $\mathsection1$]) Let $A$ be an abelian variety defined over a field $K$ finitely generated over $\mathbb{Q}$ and let $\Gamma_K=\mathrm{Gal}(\overline{K}/K)$. If $\rho_l:\Gamma_K\rightarrow \mathrm{GL}(T_l(A))$ is the $l$-adic representation of $\Gamma_K$, then the Lie algebra $\mathfrak{g}_l$ of $\rho_l(\Gamma_K)$ is algebraic and the rank of $\mathfrak{g}_l$ is independent of the prime $l$.\\

Since $V_l:=T_l(A)\otimes_{\mathbb{Z}_l}\mathbb{Q}_l$ is a semisimple $\Gamma_K$-module (Faltings and W$\ddot{\mathrm{u}}$stholz [3 Chap. 6]), the action of the Zariski closure of $\rho_l(\Gamma_K)$ in $\mathrm{GL}_{V_l}$ is also semisimple on $V_l$. Therefore it is a reductive algebraic group (Borel [1]). By Theorem 1.1, $\mathfrak{g}_l$ is algebraic. So the rank of $\mathfrak{g}_l$ is just the dimension of maximal tori. We state two more theorems:\\

\textbf{Theorem 1.2.} (Faltings and W$\ddot{\mathrm{u}}$stholz [3 Chap. 6]) Let $A$ be an abelian variety defined over a field $k$ finitely generated over $\mathbb{Q}$ and let $\Gamma_k=\mathrm{Gal}(\overline{k}/k)$. Then the map $\mathrm{End}_k(A)\otimes_{\mathbb{Z}}\mathbb{Q}_l\rightarrow \mathrm{End}_{G_k}(V_l(A))$ is an isomorphism.\\

\textbf{Theorem 1.3.} (Zarhin [6 $\mathsection5$]) Let $V$ be a finite dimensional vector space over a field of characteristic $0$. Let $\mathfrak{g}_1\subset \mathfrak{g}_2\subset \mathrm{End}(V)$ be Lie algebras of reductive subgroups of $\mathrm{GL}_V$. We assume that the centralizers of $\mathfrak{g}_1$ and $\mathfrak{g}_2$ in $\mathrm{End}(V)$ are equal and that the ranks of $\mathfrak{g}_1$ and $\mathfrak{g}_2$ are equal. Then $\mathfrak{g}_1=\mathfrak{g}_2$.\\

We are now able to prove our main theorem.\\

\textbf{Theorem 1.4.} The set $X_{\rho_{E,l}}$ is independent of $l$.\\

\textbf{Proof.} Suppose $x\in X^{cl}\backslash X_{\rho_l}$, then $(\mathfrak{g}_l)_x=\mathfrak{g}_l$. It suffices to show $\mathfrak{g}_{l'}=(\mathfrak{g}_{l'})_x:=\mathrm{Lie}(\rho_{l'}(x(\Gamma_{\textbf{k}(x)})))$ for any prime number $l'$. Since base change with finite field extension of $\textbf{k}(x)$ does not change the Lie algebras, $\mathrm{End}_{\overline{k}}(E_x)$ is finitely generated, and we have the exponential map from Lie algebras to Lie groups, we may assume that $\mathrm{End}_{\overline{k}}(E_x)=\mathrm{End}_k(E_x)$ and $\mathrm{End}_{\Gamma_k}(V_l(E_x))=\mathrm{End}_{(\mathfrak{g}_l)_x}(V_l(E_x))$. We do the same for the abelian variety $E_{\eta}/K$. We therefore have 
\[ \begin{array}{lcl}

\mathrm{dim}_{\mathbb{Q}_{l'}}(\mathrm{End}_{\mathfrak{g}_{l'}}(V_p(E_\eta)))
\stackrel{1}{=}\mathrm{dim}_{\mathbb{Q}_{l'}}(\mathrm{End}_K(E_\eta)\otimes_{\mathbb{Z}}\mathbb{Q}_{l'})\\
\\
\stackrel{2}{=}\mathrm{dim}_{\mathbb{Q}_l}(\mathrm{End}_K(E_\eta)\otimes_{\mathbb{Z}}\mathbb{Q}_l)
\stackrel{3}{=}\mathrm{dim}_{\mathbb{Q}_l}(\mathrm{End}_{\mathfrak{g}_l}(V_l(E_\eta)))\\
\\
\stackrel{4}{=}\mathrm{dim}_{\mathbb{Q}_l}(\mathrm{End}_{(\mathfrak{g}_l)_x}(V_l(E_x)))
\stackrel{5}{=}\mathrm{dim}_{\mathbb{Q}_l}(\mathrm{End}_k(E_x)\otimes_{\mathbb{Z}}\mathbb{Q}_l)\\
\\
\stackrel{6}{=}\mathrm{dim}_{\mathbb{Q}_{l'}}(\mathrm{End}_k(E_x)\otimes_{\mathbb{Z}}\mathbb{Q}_p)
\stackrel{7}{=}\mathrm{dim}_{\mathbb{Q}_{l'}}(\mathrm{End}_{(\mathfrak{g}_{l'})_x}(V_{l'}(E_x))).

\end{array}
\]
Theorem 1.2 implies the first, third, fifth and seventh equality. The dimensions of $\mathrm{End}_K(E_\eta)\otimes_{\mathbb{Z}}\mathbb{Q}_l$ and $\mathrm{End}_k(E_x)\otimes_{\mathbb{Z}}\mathbb{Q}_l$ as vector spaces are independent of $l$ imply the second and the sixth equality. $\mathfrak{g}_l=(\mathfrak{g}_l)_x$ implies the fourth equality.\\

We have $\mathrm{End}_{\mathfrak{g}_{l'}}(V_{l'}(E_\eta))=\mathrm{End}_{(\mathfrak{g}_{l'})_x}(V_{l'}(E_x))$ because the left one is contained in the right one. In other words, the centralizer of $(\mathfrak{g}_{l'})_x$ is equal to the centralizer of $\mathfrak{g}_{l'}$. We know that $(\mathfrak{g}_{l'})_x\subset \mathfrak{g}_{l'}$ are both reductive from the semisimplicity of Galois representaion (Faltings and W$\ddot{\mathrm{u}}$stholz [3 Chap. 6]). By Theorem 1.1 on $l$-independence of reductive ranks and $\mathfrak{g}_l=(\mathfrak{g}_l)_x$, we have:
\begin{equation*}
\mathrm{rank}(\mathfrak{g}_{l'})=\mathrm{rank}(\mathfrak{g}_l)=\mathrm{rank}(\mathfrak{g}_l)_x=\mathrm{rank}(\mathfrak{g}_{l'})_x.
\end{equation*}
Therefore, by Theorem 1.3 we conclude that $(\mathfrak{g}_{l'})_x= \mathfrak{g}_{l'}$ and thus prove the theorem.~$\square$\\

\textbf{Corollary 1.5 (Conjecture 5.5 [2]).} Let $k$ be a field finitely generated over $\mathbb{Q}$, $X$ a smooth, separated, geometrically connected curve over $k$ with quotient field $K$. Let $\eta$ be the generic point of $X$ and $E$ an abelian scheme over $X$. Let $\rho_l:\pi_1(X)\rightarrow \mathrm{GL}(T_l(E_{\eta}))$ be the $l$-adic representation. Then there exists a finite subset $X_E\subset X(k)$ such that for any prime $l$, $X_{\rho_{E,l}}=X_E$, where $X_{\rho_{E,l}}$ is the set of all $x\in X(k)$ such that $(G_l)_x$ is not open in $G_l:=\rho_l(\pi_1(X))$.\\

\textbf{Proof.} The uniform open image theorem for GSRP representations [2 Thm. 1.1] implies the finiteness of $X_{\rho_{E,l}}$. Theorem 1.4 implies $l$-independence.~$\square$\\

\textbf{Corollary 1.6.} Let $A$ be an abelian variety of dimension $n\geq 1$ defined over a field $K$ finitely generated over $\mathbb{Q}$. Let $\Gamma_K=\mathrm{Gal}(\overline{K}/K)$ denote the absolute Galois group of $K$. For each prime number $l$, we have the Galois representation $\rho_{l}:\Gamma_K\rightarrow \mathrm{GL}(T_l(A))$ where $T_l(A)$ is the $l$-adic Tate module of $A$. If the Mumford-Tate conjecture for abelian varieties over number fields is true, then there is an algebraic subgroup $H$ of $\mathrm{\textbf{GL}}_{2n}$ defined over $\mathbb{Q}$ such that $\rho_l({\Gamma_K})^\circ$ is open in $H(\mathbb{Q}_l)$ for all $l$.\\

\textbf{Proof.} There exists an abelian scheme $E$ over a variety $X$ defned over a number field $k$ such that the function field of $X$ is $K$ and $E_{\eta}=A$ where $\eta$ is the generic point of $X$ (see, e.g., Milne [4 $\mathsection20$]). By [5 $\mathsection1$], there exists a closed point $x\in X$ such that $(\mathfrak{g}_l)_x= \mathfrak{g}_l$. Therefore, we have $(\mathfrak{g}_l)_x= \mathfrak{g}_l$ for any prime $l$ by Theorem 1.4. Since all Lie algebras are algebraic (Theorem 1.1), if we take $H$ as the Mumford-Tate group of $E_x$, $\rho_l({\Gamma_K})^\circ$ is then open in $H(\mathbb{Q}_l)$ for all $l$.~$\square$\\

\textbf{Question.} Is the algebraic group $H$ in Corollary 1.6 isomorphic to the Mumford-Tate group of the abelian variety $A$?\\\\

\textbf{Acknowledgement.} This work grew out of an attempt to prove Conjecture 5.5 in [2] suggested by my advisor, Professor Michael Larsen. I am grateful to him for the suggestion, guidance and encouragement. I would also like to thank Anna Cadoret for her useful comments on an earlier version of this note.

\bibliographystyle{amsplain}

\bigskip
\bigskip

Department of Mathematics, Indiana University, Bloomington, IN 47405, USA\\
E-mail address: \texttt{chhui@umail.iu.edu}

\end{document}